\documentclass[runningheads]{llncs}
\usepackage{amsfonts}

\usepackage{latexsym}
\usepackage{amssymb}
\usepackage{epsfig}
\usepackage{amsmath}
\usepackage[mathscr]{eucal}
\usepackage{subfigure}
\usepackage{graphicx}

\newcommand{\ind}{1\hspace{-2.3mm}{1}}

\begin{document}

\title{Recursive methods for some problems in coding and random permutations}
\titlerunning{Recursive methods in coding}

\author{ \textbf{Ghurumuruhan Ganesan}}
\authorrunning{G. Ganesan}
\institute{Institute of Mathematical Sciences, HBNI, Chennai\\
\email{gganesan82@gmail.com }}

\date{}
\maketitle

\begin{abstract}
In this paper, we study three applications of recursion to problems in coding and random permutations. First, we consider locally recoverable codes with partial locality and use recursion to estimate the minimum distance of such codes. Next we consider weighted lattice representative codes and use recursive subadditive techniques to obtain convergence of the minimum code size. Finally, we obtain a recursive relation involving cycle moments in random permutations and as an illustration, evaluate recursions for the mean and variance.

\vspace{0.1in} \noindent \textbf{\em Key words:} Locally recoverable codes, partial locality, minimum distance, lattice identification codes, minimum size, random permutations, cycle moments.

\end{abstract}

\bigskip

\setcounter{equation}{0}
\renewcommand\theequation{\arabic{section}.\arabic{equation}}
\section{Introduction} \label{intro}
Recursive techniques are used quite frequently in coding to obtain bounds on code sizes. As a typical example, the Singleton bound~\cite{huff} obtains bounds on sizes of~\(n-\)length codes by reducing the problem to that of an~\((n-1)-\)length code. Similarly, recursive relations are also frequent in terms related to permutations like for example, Stirling numbers of the first kind~\cite{grah}. In this paper, we study further applications of recursive methods to problems in coding and random permutations. 

The paper is organized as follows: In Section~\ref{loc_rec}, we consider locally recoverable codes with partial locality and estimate the minimum distance of such codes (Theorem~\ref{main_thm}) using iterations on subcodes. Next, in Section~\ref{pf_prop1}, we study lattice representative codes with weights and prove asymptotic convergence of the minimum size, using subadditive techniques (Theorem~\ref{conv_thm}). Finally, in Section~\ref{random_perm}, we establish a recursion for cycle moments of random permutations (Theorem~\ref{cyc_mom}) and illustrate our result for the cases of mean and variance (Corollary~\ref{cor_cyc}).

\setcounter{equation}{0}
\renewcommand\theequation{\arabic{section}.\arabic{equation}}
\section{Locally recoverable codes with partial locality} \label{loc_rec}
Locally recoverable codes for erasures have tremendous applications in distributed storage and retrieval~\cite{rash} and it is therefore important to understand the properties of such codes. Typically each erasure correction is performed using a locality set of small size and it is of interest to design codes capable of correcting multiple erasures simultaneously. Such codes are also known as locally repairable codes and storage-bandwidth tradeoff and construction of such codes has been well-studied; for an overview we refer to the recent survey~\cite{balaji}. For distinction, we refer to codes above as \emph{fully} locally recoverable codes since \emph{every} symbol position has a locality set of small size associated with it. In~\cite{gop}, bounds are obtained for the minimum distance of linear fully locally recoverable codes in terms of the size of the locality sets. Later~\cite{for} studied bounds on the minimum distance of non-linear systematic fully locally recoverable codes.

In this section,  we study minimum distance of locally recoverable codes with partial locality. We assume that only a \emph{subset} of symbol positions have locality set size at most~\(r\) and obtain bounds on the minimum distance. Let~\(n\geq k \geq 1\) be integers and let~\({\cal A}\) be a set of cardinality~\(\#{\cal A} = q.\) A set~\({\cal C} \subseteq {\cal A}^{n}\) of cardinality~\(q^{k}\) is defined to be an~\((n,k)-\)code.

For a set~\({\cal U} \subseteq \{1,2,\ldots,n\}\) and an integer~\(j \in \{1,2,\ldots,n\} \setminus {\cal U},\) we say that position~\(j\) is determined by~\({\cal U}\) if there exists a function~\(g_j\) such that
\begin{equation}\label{reach_def}
c_{j} = g_{j}\left(c_i : i \in {\cal U}\right) =: g_j\left({\cal U}\right)
\end{equation}
for all codewords~\(\mathbf{c} = (c_1,\ldots,c_n) \in {\cal C}.\) In words, the symbol at the~\(j^{th}\) position of any codeword can be determined from the symbols with positions in~\({\cal U}.\) The set~\({\cal F}_{\cal U}\) of all positions determined by~\({\cal U}\) is called the \emph{reach} of~\({\cal U}.\) For integer~\(1 \leq w \leq n,\) we define
\begin{equation}\label{reach_def2}
L(w) = L(w,{\cal C}) := \max_{{\cal U} : \#{\cal U} = w} \#{\cal F}_{\cal U}.
\end{equation}
For any~\(w \geq 1\) we have that~\(L(w,{\cal C}) \leq \Delta(w),\) where~\(\Delta(w) = q^{w}\) if~\({\cal C}\) is a linear code and~\(\Delta(w) = q^{q^{w}}\) otherwise. We remark here that~\(q^{q^{w}}\) is the \emph{total} number of maps from~\({\cal A}^{w}\) to~\({\cal A}.\)

\begin{definition}
For integers~\(\tau,r \geq 2\) and a subset~\(\Theta \subseteq \{1,2,\ldots,n\},\) we say that the code~\({\cal C}\) has~\((\Theta, \tau,r)-\)local correction capability if for every subset~\({\cal P} \subseteq \Theta\) of size~\(\tau,\) there exists a set~\({\cal T}_{\cal P}  \subseteq \{1,2,\ldots,n\} \setminus {\cal P}\) of size at most~\(r\) such that each position in~\({\cal P}\) is determined by~\({\cal T}_{\cal P}.\)
\end{definition}
We define~\({\cal T}_{\cal P}\) to be the~\(r-\)locality set corresponding to the set~\({\cal P}.\) Also if~\(\Theta = \{1,2,\ldots,n\},\) we say that~\({\cal C}\) has~\((\tau,r)-\)local correction capability.

For example, consider the binary linear code~\({\cal C}\) formed in the following way: For~\(k \geq 10\) and a word~\((c_1,\ldots,c_k) \in \{0,1\}^{k},\)
let~\(d_{i_1,i_2,i_3} := c_{i_1} \oplus c_{i_2} \oplus c_{i_3}\) where~\(\oplus\) denotes addition modulo two.
There are~\({k \choose 3} =: n-k-1\) such terms~\(\{d_{i_1,i_2,i_3}\}\) which we relabel as~\(c_{k+1},\ldots,c_{n-1}.\)
Finally we let~\(c_n := \oplus_{i=1}^{k} c_i.\) We let~\((c_1,\ldots,c_n)\) be the codeword corresponding to the word~\((c_1,\ldots,c_k).\)
The collection of codewords~\({\cal C}\) has~\((\Theta,\tau,r)-\)local correction capability with~\(\Theta = \{1,2,\ldots,n-1\},\tau = 1\) and~\(r = 3.\)
For example to recover~\(c_1,\) we use the relation~\[d_{1,2,3} \oplus d_{1,2,4} \oplus d_{1,3,4} = c_1.\] In general, each bit~\(c_j, 1 \leq j \leq k\)
can be recovered in a similar manner. Because~\(k \geq 10,\) the bit~\(c_{n}\) cannot be recovered by using any three bits of~\(\{c_1,\ldots,c_{n-1}\}.\)

Let~\({\cal C}\) be any code with~\((\Theta,\tau,r)-\)local correction capability. For words~\(\mathbf{x} = (x_1,\ldots,x_n) \) and~\(\mathbf{y} = (y_1,\ldots,y_n)\) in~\({\cal C}\)
we define the Hamming distance between~\(\mathbf{x}\) and~\(\mathbf{y}\) to be~\(d(\mathbf{x},\mathbf{y}) := \sum_{i=1}^{n} \ind(x_i \neq y_i),\)
where~\(\ind(.)\) denotes the indicator function. The minimum distance of~\({\cal C}\) is then defined as~\(d({\cal C}) := \min_{\mathbf{x},\mathbf{y} \in {\cal C}} d(\mathbf{x},\mathbf{y}).\)
We have the following result.
\begin{theorem}\label{main_thm} Let~\({\cal C}\) be any~\((n,k)-\)code with~\((\Theta,\tau,r)-\)parallel correction capability and let~\(\theta = \#\Theta.\)  The minimum distance of~\({\cal C}\) satisfies
\begin{equation}\label{d_min_main}
d({\cal C}) \leq n-k+1 - T \cdot \tau,
\end{equation}
where~\(T\) is the largest integer~\(t\) such that
\begin{equation}\label{t_def_main}
t \cdot r \leq k-1+ \theta-n \text{ and } t \cdot r + \Delta(t\cdot r) \leq \theta-\tau +1.
\end{equation}
\end{theorem}
To obtain the bound~(\ref{t_def_main}), we proceed as in~\cite{for} and iteratively construct a sequence of codes with decreasing size, until no further reduction is possible. We use the pigeonhole principle at the end of each step and obtain the sufficient conditions that allow continuation of the iteration procedure.  In the proof below, we see that the first estimate in~(\ref{t_def_main}) determines the maximum number of iterations the procedure can proceed before we run out of codewords to choose from and the second estimate in~(\ref{t_def_main}) determines the maximum number of  iterations for which we are able to choose a ``fresh" locality set.

Finally, we recall that~\(n-k+1\) is the Singleton bound~\cite{huff} and is the maximum possible minimum distance of an~\((n,k)-\)code. Therefore the parameter~\(T\) is in some sense, the ``cost" for requiring partial locality.
\subsubsection*{Proof of Theorem~\ref{main_thm}}
Let~\({\cal P}_1 \subseteq \Theta\) be any set of size~\(\tau\)
and let~\({\cal I}_1 := {\cal T}_{{\cal P}_1} = \{l_1,\ldots,l_{m_1}\},m_1 \leq r\)
be the corresponding locality set of cardinality at most~\(r\) as defined in the paragraph following~(\ref{reach_def2})
that determines the value of the symbols in positions in~\({\cal P}_1.\)
For~\(\mathbf{x} = (x_1,\ldots,x_{m_1}) \in {\cal A}^{m_1},\)
let~\({\cal C}(\mathbf{x}) = {\cal C}\left(\mathbf{x},{\cal I}_1\right)\) be the set of codewords of~\({\cal C}\) such that the symbol in position~\(l_j\) equals~\(x_j\) for~\(1 \leq j \leq m_1.\)

The number of choices for~\(\mathbf{x}\) is at most~\(q^{m_1}\) and there are~\(q^{k}\) codewords in~\({\cal C}.\)
Therefore by pigeonhole principle, there exists~\(\mathbf{x}_1\) such that
\begin{equation}\label{c_one_est}
\#{\cal C}(\mathbf{x}_1) \geq \frac{\#{\cal C}}{q^{m_1}} = q^{k-m_1}.
\end{equation}
We set~\({\cal C}_1 := {\cal C}(\mathbf{x}_1)\) and let~\({\cal J}_1 := {\cal F}_{{\cal I}_1}\) be the reach of~\({\cal I}_1\) (see~(\ref{reach_def})) with cardinality~\(\tau \leq \#{\cal J}_1 \leq \Delta(r).\) The first inequality is true since~\({\cal P}_1 \subseteq {\cal J}_1\)  and the second estimate follows from~(\ref{reach_def}). By construction, all words in the code~\({\cal C}_1\) have the same values in the symbol positions determined by~\({\cal J}_1;\) i.e., if~\(a = (a_1,\ldots,a_n)\) and~\(b = (b_1,\ldots,b_n)\) both belong to~\({\cal C}_1,\) then~\(a_j = b_j\) for all~\(j \in {\cal J}_1.\)

We now repeat the above procedure with the code~\({\cal C}_1\) assuming that~\(r+n-\theta < k,\) where~\(\theta := \#\Theta.\) If~\({\cal R}_1 := \Theta \setminus \left({\cal I}_1 \cup {\cal J}_1\right)\) is the set of positions not encountered in the first iteration then
\begin{equation}\label{r1_est}
\#{\cal R}_1 \geq \#\Theta- \#{\cal J}_1 -\#{\cal I}_1 \geq \theta - \Delta(r) - r
\end{equation}
since~\(\#{\cal J}_1 \leq \Delta(r).\) For a set~\({\cal P} \subseteq {\cal R}_1\) of size~\(\tau,\) let~\({\cal I}({\cal P}) := {\cal T}_{\cal P} \bigcap {\cal R}_1\) be the union of positions within the locality sets of the selected~\(\tau\) positions in~\({\cal P},\) not encountered in the first iteration.

Suppose for every~\({\cal P} \subseteq {\cal R}_1,\) we have~\({\cal I}({\cal P}) = \emptyset.\) This means that all symbols with positions in~\({\cal R}_1\) can simply be determined by the symbol values with positions in~\({\cal I}_1 \cup {\cal J}_1.\) This in turn implies that the symbols with positions in~\({\cal I}_1\) determine \emph{all} the symbols with positions in~\(\Theta.\) Using~\(r+n-\theta < k\) we then get that the total number of words in the code~\({\cal C}\) is at most~\[q^{\#{\cal I}_1} \cdot q^{n-\#\Theta} = q^{m_1+n-\theta} \leq q^{r+n-\theta} < q^{k},\] a contradiction. Thus there exists~\({\cal P}_2 \subseteq {\cal R}_1\) of size~\(\tau\) whose corresponding set~\({\cal I}_2 := {\cal I}({\cal P}_2)\) is not completely  contained in~\({\cal I}_1 \cup {\cal J}_1.\)

Letting~\(1 \leq m_2 \leq r\) denote the cardinality of~\({\cal I}_2\)  and using the pigeonhole principle as before, we get a code~\({\cal C}_2 \subseteq {\cal C}_1\) of size~\[\#{\cal C}_2 \geq \frac{\#{\cal C}_1}{q^{m_2}} \geq q^{k-m_1-m_2} \geq q^{k-2r}\] and all of whose words have the same symbol values in the positions determined by~\({\cal I}_1 \cup {\cal I}_2.\) In the above, we use the estimate for~\({\cal C}_1\) obtained in~(\ref{c_one_est}). As before, let~\({\cal J}_2 \subseteq \{1,2,\ldots,n\}\) be the set of positions of the codeword symbols determined by the set~\({\cal I}_1 \cup {\cal I}_2 \) so that~\({\cal P}_1 \cup {\cal P}_2 \subseteq {\cal J}_2.\)
The set~\({\cal I}_1 \cup {\cal I}_2\) has cardinality at most~\(2r\) and so we have from~(\ref{reach_def2}) that the reach~\({\cal J}_2\) has cardinality~\[2\tau \leq \#{\cal P}_1 + \#{\cal P}_2 \leq \#{\cal J}_2 \leq \Delta(2r).\]

If~\({\cal R}_2 := \Theta \setminus \left({\cal I}_1 \cup {\cal I}_2 \cup {\cal J}_2\right),\) then~\(\#{\cal R}_2 \geq  \theta - 2r  - \Delta(2r).\) Continuing this way, after the end of~\(t\) iterations, we have a code~\({\cal C}_t\) of size
\begin{equation}\label{ct_est}
\#{\cal C}_t \geq q^{k - \sum_{j=1}^{t} \#{\cal I}_j} \geq q^{k - t \cdot r}
\end{equation}
and a set~\({\cal R}_t \subseteq \Theta\) of remaining positions not fixed so far, with cardinality
\begin{equation}\label{rt_est}
\#{\cal R}_t  \geq \theta - \sum_{j=1}^{t} \#{\cal I}_j - \#{\cal J}_t \geq \theta - t\cdot r - \Delta(t\cdot r) .
\end{equation}
The above procedure can therefore be performed for at least~\(T\) steps where~\(T\) is the largest integer~\(t\) such that
\begin{equation}\label{t_def2}
t \cdot r \leq k-1+\theta-n \text{ and } t\cdot r + \Delta(t\cdot r) \leq \theta-\tau+1.
\end{equation}
The first condition in~(\ref{t_def2}) ensures that~\(k - T \cdot r \geq 1\) and so the code~\({\cal C}_T\) has at least two codewords.
The second condition in~(\ref{t_def2}) ensures that the set~\({\cal P}_j \subseteq \Theta\) of symbols we pick is at least~\(\tau\) and so
\begin{equation}\label{jt_est}
\#{\cal J}_j \geq \sum_{l=1}^{j} \#{\cal P}_l \geq j \cdot\tau
\end{equation}
for each~\(1 \leq j \leq T.\)


Since~\({\cal C}_T \subseteq {\cal C},\) the minimum distance~\(d({\cal C}_T)\) of~\({\cal C}_T\) is at least the minimum distance~\(d({\cal C})\) of~\({\cal C}.\) By definition we recall that the symbol values in positions determined by the set~\(\bigcup_{1 \leq j \leq T} {\cal I}_j \bigcup {\cal J}_T := \{1,2,\ldots,n\} \setminus {\cal Q}_T\) is the same for all the words in~\({\cal C}_T.\) For every word~\(\mathbf{x}  = (x_1,\ldots,x_n) \in {\cal C}_T,\) we therefore let~\(\mathbf{x}_T = (x_i)_{i \in {\cal Q}_T}\) be the reduced word obtained by just considering the symbols in the remaining positions determined by~\({\cal Q}_T.\) Defining the reduced code~\({\cal D}_T = \{\mathbf{x}_T : \mathbf{x} \in {\cal C}_T\}\) we then have that the minimum distance~\(d({\cal D}_T) \geq d({\cal C}_T) \geq d({\cal C}).\)

The length of the each word in~\({\cal D}_T\) equals~\(n - \#{\cal Q}_T\) and so using the estimate for~\(\#{\cal D}_T = \#{\cal C}_T\) from~(\ref{ct_est}) and the Singleton bound we have
\begin{equation}
d({\cal D}_T) \leq \left(n - \sum_{j=1}^{T} \#{\cal I}_j - \#{\cal J}_T\right)  -  \left(k - \sum_{j=1}^{T} \#{\cal I}_j\right) + 1. \nonumber
\end{equation}
Thus~\(d({\cal D}_T) \leq n - k+1 - \#{\cal J}_T \leq n - k+1 - T \cdot \tau,\) by~(\ref{jt_est}). Using the fact that~\(d({\cal C}) \leq d({\cal D}_T),\) we then get~(\ref{d_min_main}).\;\;\;\;\;\;\;\;\;\;\;\;\;\;\;\;\;\;\;\;\;\;\;\;\;\;\;\;\;\;\;\;\;\;\;\;\;\;\;\;\;\;\;\;\;\;\;\;\;\;\;\;
\;\;\;\;\;\;\;\;\;~\(\qed\)


\setcounter{equation}{0}
\renewcommand\theequation{\arabic{section}.\arabic{equation}}
\section{Lattice representative codes} \label{pf_prop1}
Representative codes~\cite{karp} (also known as hitting sets in some contexts) are important from both theoretical and application perspectives. In~\cite{lin} the minimum size of hitting sets that intersect all combinatorial rectangles of a given volume are studied. Explicit constructions were described using expander graphs and random walks. Later~\cite{chan} determined lower bounds for the hitting set size of combinatorial rectangles and also illustrated an application in approximation algorithms. Recently~\cite{bhas} used fractional perfect hash families to study construction of explicit hitting sets for combinatorial shapes.


In this section, we study lattice representative codes for \emph{weighted} rectangles where each vertex is assigned a positive finite weight. We study the minimum size of a representative code that intersects \emph{all} subsets of a given minimum weight. For integers~\(d, m \geq 1\) let~\({\cal S}  = {\cal S}(m):= \{1,2,\ldots,m\}^{d}.\) We refer to elements of~\({\cal S}\) as points and for a point~\(v = (v_1,\ldots,v_d) \in {\cal S},\) we refer to~\(v_i\) as the~\(i^{th}\) entry. For each point~\(v \in {\cal S},\) we assign a finite positive weight~\(w(v).\)  For a set~\({\cal U} \subseteq {\cal S}\) the corresponding weight~\(w\left({\cal U}\right) := \sum_{v \in {\cal U}} w(v)\) is the sum of the weights of the points in~\({\cal U}.\) The size of~\({\cal U}\) is the number of points in~\({\cal U}\) and is denoted by~\(\#{\cal U}.\)

Let
\begin{equation}\label{w_low_def}
1= \inf_{m} \min_{v \in {\cal S}(m)} w(v) \leq \sup_m \max_{v \in {\cal S}(m)} w(v)  =: \beta < \infty
\end{equation}
and for~\(0 < \epsilon <1\) say that~\({\cal B} \subseteq {\cal S}\) is an \(\epsilon-\)representative code or simply representative code if~\({\cal B} \cap {\cal U} \neq \emptyset\) for any set~\({\cal U} \subseteq {\cal S}\) of weight~\(w\left({\cal U}\right) \geq \epsilon~\cdot~m^{d}.\)

The following result obtains an estimate on the minimum size~\(b_m\) of an\\\(\epsilon-\)representative code.
\begin{theorem}\label{conv_thm}
For any~\(0 < \epsilon <1\) and~\(\beta \geq 1\) we have that
\begin{equation}\label{bm_bounds}
m^{d}\cdot (1-\epsilon) \leq b_m \leq m^{d}\cdot \left(1-\frac{\epsilon}{\beta}\right)+1.
\end{equation}

Suppose the weight function satisfies the following monotonicity relation: If~\(u\) and~\(v\) are any two points of~\({\cal S}\) differing only in the~\(i^{th}\) entry and~\(u_i > v_i,\) then the weights~\(w(u) \leq w(v).\) We then have
\begin{equation}\label{b_conv}
\frac{b_m}{m^{d}} \longrightarrow \lambda
\end{equation}
as~\(m \rightarrow \infty\) where~\(1-\epsilon \leq \lambda  \leq 1-\frac{\epsilon}{\beta}.\)
\end{theorem}
Thus there exists a \emph{fraction} of vertices in a large rectangle that hits \emph{all} sets of a given minimum weight. Moreover, if the weight assignment is monotonic, then the scaled minimum representative code size converges to a positive constant strictly between~\(0\) and~\(1.\)

An example of  non-trivial weight assignment that satisfies the monotonicity relation is the following: Defining~\(w(1,1) := 2\) we iteratively assign the weight of each vertex in the set~\(\{1,2,\ldots,i+1\}^{d} \setminus \{1,2,\ldots,i\}^d\) as~\(1+\frac{1}{i}.\) The conditions in Theorem~\ref{conv_thm} are then satisfied with~\(\beta =2.\)

\subsubsection*{Proof of Theorem~\ref{conv_thm}}
We begin with the proof of~(\ref{bm_bounds}). Throughout we assume that~\(d =2\) and an analogous analysis holds for general~\(d.\) If~\({\cal F}\) is any representative code of~\({\cal S},\) then by definition, the weight of the set~\({\cal S} \setminus {\cal F}\) is at most~\(\epsilon m^{2}\) and  since the weight of each vertex is at least one (see~(\ref{w_low_def})), we get that the number of points in~\({\cal S} \setminus {\cal F}\) is at most~\(\epsilon m^{2}.\) This implies that the size of~\({\cal F}\) is at least~\((1-\epsilon) m^2\) and so~\(b_m \geq (1-\epsilon) m^{2}.\)

To find an upper bound on~\(b_m,\) we let~\({\cal T} \subseteq {\cal S}\) be any ``critical" set such that the weight of~\({\cal T}\) is at most~\(\epsilon m^{2}-1\) and the weight of~\({\cal T} \cup \{v\}\) for any point~\(v \in {\cal S} \setminus {\cal T}\) is at least~\(\epsilon m^{2}.\)
The set~\({\cal S} \setminus {\cal T}\)  is then a representative code of~\({\cal S}\) and since the weight of any point is at most~\(\beta\) (see~(\ref{w_low_def})), we get that the number of vertices in~\({\cal T}\) is at least~\(\frac{\epsilon }{\beta} \cdot m^{2} - 1.\)  This in turn implies that~\(b_m \leq m^{2}\left(1-\frac{\epsilon}{\beta}\right)+1.\) This proves~(\ref{bm_bounds}).

To prove~(\ref{b_conv}), we use a subsequence argument analogous to the proof of Fekete's lemma. For integers~\(m \geq r \geq 1,\) we let~\(m = k \cdot r + s,\) where~\(k\geq 1\) and~\(0 \leq s \leq r-1\) are integers and split~\(\{1,2,\ldots,m\}^2\) into four sets
\[{\cal S}_1:=\{1,2,\ldots,kr\}^2,\;\;\;{\cal S}_2 := \{1,2,\ldots,kr\} \times \{kr+1,kr+2,\ldots,kr+s\},\]
\[{\cal S}_3 := \{kr+1,\ldots,kr+s\} \times \{1,2,\ldots,kr\} \text{ and } {\cal S}_4 := \{kr+1,\ldots,kr+s\}^2.\] Thus~\({\cal S}_2\) is essentially a ``rotated" version of~\({\cal S}_3.\) For~\(1 \leq i \leq 4,\) let~\({\cal G}_i(\epsilon)\) be a representative code of~\({\cal S}_i\) and let~\({\cal R}\) be any set in~\(\{1,2,\ldots,m\}^2\) of weight~\(w\left({\cal R}\right) \geq \epsilon m^2 = \epsilon(kr+s)^2.\) We first see that~\(\bigcup_{i=1}^{4} {\cal G}_i(\epsilon)\) is a representative code of~\(\{1,2,\ldots,m\}^2.\) Indeed if~\({\cal R}_i = {\cal R} \cap {\cal S}_i,\) then using the fact that~\(w\left({\cal R}\right) = \sum_{i=1}^{4} w\left({\cal R}_i\right),\) we get that either~\(w\left({\cal R}_1\right) \geq \epsilon (kr)^2\) or~\(w\left({\cal R}_2\right) \geq \epsilon k r s\) or~\(w\left({\cal R}_3\right) \geq \epsilon k r s\) or~\(w\left({\cal R}_4\right) \geq \epsilon s^2.\) Consequently, we must have that~\({\cal R} \bigcap \bigcup_{i=1}^{4} {\cal G}_i(\epsilon) \neq \emptyset.\)

If~\(b^{(i)}\) denotes the minimum size of a representative code of~\({\cal S}_i,\) then from the discussion above we get
\begin{equation}\label{bkr}
b_{kr+s} \leq \sum_{i=1}^{4}b^{(i)} \leq b^{(1)} + 2krs + s^2 \leq b^{(1)} + (2k+1)r^2
\end{equation}
where the second inequality in~(\ref{bkr}) follows from the trivial estimate that the size of any representative code of~\({\cal R}_i\) is at most the total number of points in~\({\cal R}_i\) and the final inequality in~(\ref{bkr}) follows from the fact that~\(s \leq r.\)

To estimate~\(b^{(1)},\) we split~\({\cal S}_1 := \{1,2,\ldots,kr\}^2\) into~\(k^2\) disjoint rectangles~\({\cal T}_i,\)\\\( 1 \leq i \leq k^2\) each containing~\(r^2\) points with~\({\cal T}_1 = \{1,2,\ldots,r\}^2.\) If~\(c_i\) denotes the minimum size of a representative code of~\({\cal T}_i,\) then using the weight monotonicity relation, we get that~\(c_i \leq c_1.\) To see this is true suppose~\({\cal T}_2 = \{1,2,\ldots,r\} \times \)\\\( \{r+1,\ldots,r+2r\}\) so that~\({\cal T}_2 = {\cal T}_1 + (r,0)\) is obtained by translation of~\({\cal T}_1.\) If~\({\cal U} \subset {\cal T}_2\) is any set of weight at least~\(\epsilon r^{2}\) then~\({\cal U} - (r,0) \subseteq {\cal T}_1\) also has weight at least~\(\epsilon r^2,\) by the weight monotonicity relation. Consequently if~\({\cal W}_1\) is a representative code of~\({\cal T}_1,\) then~\({\cal W}_1 + (r,0)\) is a representative code of~\({\cal T}_2.\) Thus~\(c_2 \leq c_1\) and the proof of general~\(c_i\) is analogous.

From~(\ref{bkr}) and the discussion in the above paragraph, we get that\\\(b_m = b_{kr+s} \leq k^2 b_r + (2k+1)r^2\) and so
\[\frac{b_m}{m^2} \leq \frac{k^2 b_r + (r^2(2k+1))}{m^2} = \left(\frac{kr}{kr+s}\right)^2\left( \frac{b_r}{r^2} + \frac{2k+1}{k^2}\right).\]
If~\(m \rightarrow \infty\) with~\(r\) fixed, then~\(k = k(m,r) = \frac{m-s}{r} \geq \frac{m-r}{r} \rightarrow \infty\) as well and so~\(\frac{kr}{kr+s} \rightarrow 1.\) This in turn implies that
\begin{equation}\label{bm_c}
\limsup_m \frac{b_m}{m^2} \leq \limsup_m \left( \frac{b_r}{r^2} + \frac{2k+1}{k^2}\right) = \frac{b_r}{r^2}.
\end{equation}
Since~\(r \geq 1\) is arbitrary we get from~(\ref{bm_c}) that~\[\limsup_m \frac{b_m}{m^2} = \liminf_{r} \frac{b_r}{r^2} = \inf_{r} \frac{b_r}{r^2} =: \lambda.\] Also, the bounds for~\(\lambda\) follow from~(\ref{bm_bounds}).
\;\;\;\;\;\;\;\;\;\;\;\;\;\;\;\;\;\;\;\;\;\;\;\;\;\;\;\;\;\;\;\;\;\;\;\;\;\;\;\;\;\;\;\;\;\;\;\;\;\;\;\;\;\;\;\;\;~\(\qed\)

\renewcommand{\theequation}{\arabic{section}.\arabic{equation}}
\setcounter{equation}{0}
\section{Random permutations} \label{random_perm}
Random permutations and applications are frequently encountered in computing problems and it is of interest to study the cycle properties of a randomly chosen permutation. The papers~\cite{gon},~\cite{shep} studied limiting distributions for the convergence of the number of cycles and cycles lengths of a uniform random permutation, after suitable renormalization. Later~\cite{arr} used Poisson approximation and estimates on the total variation distance to study the convergence of the overall cycle structure to a process of independent Poisson random variables. Recently~\cite{betz} have used probability generating functions to study convergence of number of cycles of uniform random permutations conditioned not to have large cycles, scaled and centred, to the Gaussian distribution.

From the combinatorial aspect, Stirling numbers of the first kind and generating functions have been used to study random permutation statistics. Using the Flajolet-Sedgewick theorem it is possible to enumerate permutations with constraints~\cite{grah}. In this section, we use conditioning to obtain a recursive relation involving cycle moments of random permutations. As an illustration, we compute recursive relation involving the mean and the variance of the number of cycles in a uniformly random permutation.

We begin with a couple of definitions. A permutation~\(\pi\) of~\(\{1,2,\ldots,n\}\) is a bijective map~\(\pi : \{1,2,\ldots,n\} \rightarrow \{1,2,\ldots,n\}.\) The total number of possible permutations of~\(\{1,2,\ldots,n\}\) is therefore~\[n! := n\cdot (n-1) \cdots 2\cdot 1.\] A \emph{cycle} of length~\(k\) in a permutation~\(\pi\) is a~\(k-\)tuple~\((i_1,\ldots,i_{k})\) such that~\(\pi(i_j) = i_{j+1}\) for~\(1 \leq j \leq k-1\) and~\(\pi(i_k) = i_1.\) Every number in~\(\{1,2,\ldots,n\}\) belongs to some cycle of~\(\pi\) and this provides an alternate representation of~\(\pi;\) for example~\((1345)(267)(89)\) is the cycle representation of the permutation~\(\pi\) on~\(\{1,2,\ldots,9\}\) satisfying~\(\pi(1) = 3, \pi(3) = 4, \pi(4) = 5, \pi(5) = 1,\pi(2) = 6, \pi(6) = 7, \pi(7) = 2, \pi(8) = 9, \pi(9) = 8.\)

Let~\(\Pi\) denote a uniformly chosen random permutation of~\(\{1,2,\ldots,n\}\) defined on the probability space~\((\Omega_n, {\cal F}_n, \mathbb{P}_n)\) so that \[\mathbb{P}_n(\Pi = \pi) = \frac{1}{n!}\] for any deterministic permutation~\(\pi.\)  Let~\(N_n = N_n(\Pi)\) be the random number of cycles in~\(\Pi\) and for integers~\(n,s \geq 1,\) set~\(\mu_{0,s} := 0\) and~\(\mu_{n,s} := \mathbb{E}N_n^{s}.\) We have the following result.
\begin{theorem}\label{cyc_mom} For integers~\(n,s \geq 1\) we have
\begin{equation}\label{rec_eq}
\mu_{n,s} = 1 + \frac{1}{n} \sum_{r=1}^{s} \sum_{j=1}^{n-1} {s \choose r} \mu_{j,r},
\end{equation}
where~\({s \choose r} = \frac{s!}{r!(s-r)!}\) is the Binomial coefficient.
\end{theorem}
From the recursive structure of equation~(\ref{rec_eq}), we then have that~\(\mu_{n,s}\) could be computed using the previous values~\(\{\mu_{j,r}\}_{j \leq n-1, r \leq s}.\)

As a Corollary of Theorem~\ref{cyc_mom} we have the following recursive relations for the mean and variance of~\(N_n.\)
\begin{corollary}\label{cor_cyc} The mean~\(\mu_n := \mu_{n,1}\) satisfies~\(\mu_1 = 1\) and the recursive equation
\begin{equation}\label{mean_cyc}
\mu_n  = 1+\frac{1}{n} \sum_{i=1}^{n-1} \mu_i
\end{equation}
for~\(n \geq 2.\) The sequence~\(H_n := \sum_{j=1}^{n}\frac{1}{j}\) is the unique sequence satisfying~(\ref{mean_cyc}).

The variance~\(v_n = var(N_n) := \mu_{n,2}-\mu_{n,1}^2\) satisfies~\(v_1 = 0\) and the recursive equation
\begin{equation}\label{var_cyc}
v_n = 1+\frac{1}{n}\sum_{i=1}^{n-1} v_i -\frac{H_n}{n}.
\end{equation}
The sequence~\(M_n := H_n - \sum_{i=1}^{n}\frac{1}{i^2}\) is the unique sequence satisfying~(\ref{var_cyc}).
\end{corollary}
Using~(\ref{mean_cyc}),~(\ref{var_cyc}) and the recursive relation~(\ref{rec_eq}), we could similarly compute higher order moments.

We prove Theorem~\ref{cyc_mom} and Corollary~\ref{cor_cyc} in that order.
\subsubsection*{Proof of Theorem~\ref{cyc_mom}} 
To obtain the desired recursive relation, we condition on the length of the first cycle and study the number of  cycles in the remaining set of elements.

Let~\({\cal S}_1\) denote the cycle of the random permutation~\(\Pi\) containing the number~\(1\) and let~\(L_1 = \#{\cal S}_1\) be the length of~\({\cal S}_1\) so that~\({\cal S}_1\) is an~\(L_1-\)tuple. If~\(L_1 = k \leq n-1,\) then~\(\Pi\) induces a permutation~\(\sigma : \{1,2,\ldots,n-k\} \rightarrow \{1,2,\ldots,n-k\}\) on the remaining~\(n-k\) numbers~\(\{1,2,\ldots,n\}\setminus {\cal S}_1\) in the following way. Arrange the numbers in~\(\{1,2,\ldots,n\} \setminus {\cal S}_1\)
in increasing order~\(j_1 < j_2 < \ldots < j_{n-k}\) and suppose that~\(\pi(j_l) = m_l\) for~\(1 \leq l \leq n-k.\) The induced permutation~\(\sigma\) then satisfies~\(m_l = j_{\sigma(l)}\) for~\(1 \leq l \leq n-k.\)

Conditional on~\(L_1 = k\) we now see that~\(\sigma\) is uniformly distributed in the sense that for any deterministic permutation~\(\sigma_0 : \{1,2,\ldots,n-k\} \rightarrow \{1,2,\ldots,n-k\}\) we have
\begin{equation}\label{cond_eq}
\mathbb{P}_n\left(\sigma = \sigma_0 | L_1 = k\right) = \mathbb{P}_{n-k}(\sigma_0) = \frac{1}{(n-k)!}.
\end{equation}
To see~(\ref{cond_eq}) is true, we first write
\begin{equation}\label{cond_exp}
\mathbb{P}_n\left(\sigma = \sigma_0 | L_1 = k\right) = \frac{\mathbb{P}_n \left(\{\sigma = \sigma_0\} \cap \{L_1 = k\}\right)}{\mathbb{P}_n(L_1=  k)}.
\end{equation}
If~\(k=1,\) then the numerator in the right side of~(\ref{cond_exp}) is~\(\frac{1}{n!}.\) Moreover, if the first cycle simply consists of the single element~\(1,\) then the remaining~\(n-1\) numbers can be arranged in~\((n-1)!\) ways and so~\(\mathbb{P}_n(L_1 = 1) = \frac{(n-1)!}{n!}.\) Thus~(\ref{cond_eq}) is true for~\(k=1.\)

For~\(2 \leq k \leq n-1,\) we have from~(\ref{cond_exp}) that~\(\mathbb{P}_n\left(\sigma = \sigma_0 | L_1 = k\right)\) equals
\begin{equation}
\frac{\sum_{(i_1,\ldots,i_{k-1})}\mathbb{P}_n \left(\{\sigma = \sigma_0 \}\cap \{{\cal S}_1 = (1,i_1,\ldots,i_{k-1})\}\right)}{\sum_{(i_1,\ldots,i_{k-1})}\mathbb{P}_n({\cal S}_1 = (1,i_1,\ldots,i_{k-1}))} \nonumber\\
\label{grt_one}
\end{equation}
where the summation is over all~\(k-1\) tuples~\((i_1,\ldots,i_{k-1})\) containing distinct elements. For any~\((1,i_1,\ldots,i_{k-1}),\) the term
\begin{equation}\label{eq_a}
\mathbb{P}_n \left(\{\sigma = \sigma_0\} \cap \{{\cal S}_1 = (1,i_1,\ldots,i_{k-1})\}\right) = \frac{1}{n!}
\end{equation}
and
\begin{equation}\label{eq_b}
\mathbb{P}_n({\cal S}_1 = (1,i_1,\ldots,i_{k-1})) = \frac{(n-k)!}{n!}
\end{equation}
since there are~\((n-k)!\) ways to permute the remaining~\(n-k\) elements of the set\\\(\{2,\ldots,n\} \setminus \{i_1,\ldots,i_{k-1}\}.\) Substituting~(\ref{eq_a}) and~(\ref{eq_b}) into~(\ref{grt_one}), we get~(\ref{cond_eq}).

Summing~(\ref{eq_b}) over all~\(k-1\) tuples with distinct entries (for which there are~\((n-1)\cdot (n-2) \cdots (n-k+1)\) choices), we also get  that~\(\mathbb{P}_n(L_1 = \#{\cal S}_1 = k) = \frac{1}{n}.\)
From the discussion in the previous paragraph, we get that the above relation holds for all~\(1 \leq k \leq n.\) Thus
\begin{equation}
\mu_{n,s} = \mathbb{E}_nN_n^{s} = \sum_{k=1}^{n} \mathbb{E}_n(N_n^s | L_1=k) \mathbb{P}_n(L_1=k)= \frac{1}{n} \sum_{k=1}^{n} \mathbb{E}_n (N_n^s | L_1=k). \label{egr2}
\end{equation}
If~\(k = n\) then~\(N_n=1\) and if~\(1 \leq k \leq n-1,\) then~\(N_n = 1+M_n,\) where~\(M_n\) is the number of cycles in the induced permutation~\(\sigma.\)
Therefore we get from~(\ref{egr2}) that
\begin{equation}\label{egr3}
\mu_{n,s} = \frac{1}{n} + \frac{1}{n}\sum_{k=1}^{n-1}\mathbb{E}_n\left((1+M_n)^s|L_1=k\right).
\end{equation}

Using the conditional distribution equivalence~(\ref{cond_eq}), we have for~\(1 \leq k \leq n-1\) that~\(\mathbb{E}_n\left((1+M_n)^s|L_1=k\right) \) equals
\begin{equation}
\mathbb{E}_{n-k}(1+N_{n-k})^s = 1+\sum_{r=1}^{s} {s \choose r} \mathbb{E}_{n-k}N^r_{n-k},
\label{egr_xx}
\end{equation}
by the Binomial expansion. Substituting~(\ref{egr_xx}) into~(\ref{egr3}) we get~(\ref{rec_eq}).\;\;\;\;\;\;\;\;\;\;~\(\qed\)

\subsubsection*{Proof of Corollary~\ref{cor_cyc}} 
We begin with the proof of~(\ref{mean_cyc}). Setting~\(s=1\) in~(\ref{rec_eq}) and~\(\mu_n = \mu_{n,1},\) we get that~\(\mu_1 = 1\) and for~\(n \geq 2,\)
we get that~\(\mu_n\) satisfies~(\ref{mean_cyc}). We first see by induction that~\(H_n\) as defined in Proposition~\ref{cor_cyc} satisfies~(\ref{mean_cyc}).
For~\(n=2,\) this statement is true and suppose~\(H_l\) satisfies~(\ref{mean_cyc}) for~\(1 \leq l \leq n-1.\)
For~\(l=n,\) the right side of~(\ref{mean_cyc}) evaluated with~\(\mu_i = H_i\) equals
\begin{equation}\label{ght}
1+\frac{1}{n} \sum_{i=1}^{n-1} \sum_{j=1}^{i}\frac{1}{j}= 1+ \frac{1}{n}\sum_{j=1}^{n-1} \sum_{i=j}^{n-1}\frac{1}{j} = 1+\frac{1}{n}\sum_{j=1}^{n-1}\frac{n-j}{j},
\end{equation}
by interchanging the order of summation in the second equality. The final term in~(\ref{ght}) equals~\(H_n\) and this proves the induction step.

Suppose now that~\(\{b_n\}\) is some sequence satisfying~(\ref{mean_cyc}) with~\(b_1 =1\) and let~\(u_n = b_n-H_n\) denote the difference.
The sequence~\(\{u_n\}\) satisfies~\(u_1 = 0\) and~\( u_n = \frac{1}{n} \sum_{i=1}^{n-1} u_i\)
for all~\(n \geq 2.\) Thus~\(u_2 = \frac{u_1}{2} = 0\) and iteratively, we get~\(u_n~=~u_2~=~0\) for all~\(n \geq 2.\) Thus~\(H_n\) is the unique sequence satisfying~(\ref{mean_cyc}).

We now obtain the variance estimate as follows. Letting~\(d_n := \mu_{n,2}\) and~\(\mu_n := \mu_{n,1} = H_n,\) we
get from~(\ref{rec_eq}) that
\begin{equation}\label{grt2}
d_n = 1+ \frac{1}{n} \sum_{i=1}^{n-1} (d_i+2\mu_i) = 2\mu_n-1 + \frac{1}{n} \sum_{i=1}^{n-1} d_i,
\end{equation}
since~\(\frac{1}{n} \sum_{i=1}^{n-1} \mu_i = \mu_n-1\) (see~(\ref{mean_cyc})). From~(\ref{grt2}) we get that~\(v_n = d_n-\mu_n^2\) equals
\begin{eqnarray}
v_n &=&\frac{1}{n}\sum_{i=1}^{n-1} (d_i-\mu_i^2) + \frac{1}{n} \sum_{i=1}^{n-1}\mu_i^2 - (\mu_n-1)^2 \nonumber\\
    &=&\frac{1}{n}\sum_{i=1}^{n-1} v_i + \frac{1}{n} \sum_{i=1}^{n-1}\mu_i^2 - (\mu_n-1)^2. \nonumber
\end{eqnarray}

It only remains to see that~\(\frac{1}{n}\sum_{i=1}^{n-1} \mu_i^2 - (\mu_n-1)^2 = 1-\frac{H_n}{n}\) and for that
we use~\(\mu_i = H_i = \sum_{j=1}^{i} \frac{1}{j}\) (see~(\ref{mean_cyc})) to first get that~\(\frac{1}{n} \sum_{i=1}^{n-1} \mu_i^2\) equals
\begin{eqnarray}
&&\frac{1}{n}\sum_{i=1}^{n-1}\sum_{j_1=1}^{i} \sum_{j_2=1}^{i} \frac{1}{j_1\cdot j_2} = \frac{1}{n} \sum_{j_1=1}^{n-1} \sum_{j_2=1}^{n-1} \sum_{i = \max(j_1,j_2)}^{n-1} \frac{1}{j_1\cdot j_2}\nonumber\\
&&\frac{1}{n} \sum_{j_1=1}^{n-1} \sum_{j_2=1}^{n-1} \frac{(n-\max(j_1,j_2))}{j_1\cdot j_2} = \frac{1}{n} \sum_{j_1=1}^{n-1} \Delta(j_1), \nonumber\\ \label{step_one}
\end{eqnarray}
where~\(\Delta(j_1) = \sum_{j_2=1}^{j_1} \frac{n-j_1}{j_1\cdot j_2} + \sum_{j_2=j_1+1}^{n-1} \frac{n-j_2}{j_1\cdot j_2}\) equals
\[\sum_{j_2=1}^{n-1}\frac{n}{j_1\cdot j_2} -\sum_{j_2=1}^{j_1} \frac{1}{j_2} - \sum_{j_2=j_1+1}^{n-1} \frac{1}{j_1}.\]
Thus~\(\frac{1}{n} \sum_{i=1}^{n-1} \mu_i^2\) equals
\begin{equation}
\sum_{j_1=1}^{n-1}\sum_{j_2=1}^{n-1} \frac{1}{j_1 \cdot j_2} - \frac{1}{n} \sum_{j_1=1}^{n-1}\sum_{j_2=1}^{j_1}\frac{1}{j_2} - \frac{1}{n}\sum_{j_1=1}^{n-1} \sum_{j_2=j_1+1}^{n-1}\frac{1}{j_1}. \label{step_two}
\end{equation}

The first term in~(\ref{step_two}) is~\[\sum_{j_1=1}^{n-1}\sum_{j_2=1}^{n-1} \frac{1}{j_1 \cdot j_2} = H_{n-1}^2 = \left(H_n-\frac{1}{n}\right)^2\]
and the second term in~(\ref{step_two}) is~\[\frac{1}{n} \sum_{j_1=1}^{n-1}\sum_{j_2=1}^{j_1}\frac{1}{j_2}  = \frac{1}{n}\sum_{j_1=1}^{n-1}H_{j_1} = H_n-1\]
using the fact that~\(\mu_n = H_n\) satisfies~(\ref{mean_cyc}). The third term  in~(\ref{step_two})
equals~\[\frac{1}{n}\sum_{j_1=1}^{n-1}\frac{n-1-j_1}{j_1}  = \left(\frac{n-1}{n}\right)\left(H_{n} - 1-\frac{1}{n}\right)\]
after rearrangement of  terms. Substituting these three expressions into~(\ref{step_two}), we get that~\(\frac{1}{n} \sum_{i=1}^{n-1} \mu_i^2 \) equals~\(1+ (H_n-1)^2 - \frac{H_n}{n},\) which is what we wanted to prove. Finally, arguing as before, we also have that~\(M_n\) is the unique sequence satisfying~(\ref{var_cyc}).\;\;\;\;\;\;\;\;\;\;\;\;\;\;\;\;\;\;\;\;\;\;\;\;\;\;\;\;\;\;\;\;\;\;\;\;\;\;\;
\;\;\;\;\;\;\;\;\;\;\;\;\;\;\;\;\;\;\;\;\;\;\;\;\;\;\;\;\;\;\;\;\;\;\;\;\;\;\;
\;\;\;\;\;\;\;\;\;\;\;\;\;\;\;\;~\(\qed\)

\subsubsection*{Acknowledgements}
I thank Professors Rahul Roy, V. Guruswami, C. R. Subramanian and the referees for crucial comments that led to an improvement of the paper.
I also thank IMSc for my fellowships.

\bibliographystyle{plain}

\end{document}